\documentclass[final,1p,times]{elsarticle}
\usepackage{txfonts}

\usepackage{amssymb}
 \usepackage{amsthm}
\usepackage{amscd}
\usepackage{amsmath}
\usepackage{amsfonts}
\usepackage{amssymb}
\usepackage{graphicx}

\numberwithin{equation}{section}

\newcommand{\ra}{\rightarrow}

\newcommand{\f}{\frac}

\newcommand{\be}{\begin{equation}}
\renewcommand{\ra}{\rightarrow}
\newcommand{\ee}{\end{equation}}
\newcommand{\bea}{\begin{eqnarray}}
\newcommand{\eea}{\end{eqnarray}}
\newcommand{\bna}{\begin{eqnarray*}}
\newcommand{\ena}{\end{eqnarray*}}

\renewcommand{\le}{\left}
\newcommand{\ri}{\right}

\newcommand{\ve}{\varepsilon}
\newcommand{\ep}{\epsilon}

\journal{$\ast\ast\ast$}

\begin{document}

\begin{frontmatter}

\title{Min-max level estimate for a singular quasilinear polyharmonic equation in $\mathbb{R}^{2m}$}

\author[label1,label2]{Liang Zhao}
 \ead[label1]{liangzhao@bnu.edu.cn}
\author[label1]{Yuanyuan Chang}
\address[label1]{School of Mathematical Sciences,
 Beijing Normal
University, Beijing 100875, P. R. China}
\address[label2]{Key Laboratory of Mathematics and
Complex Systems of Ministry of Education, Beijing Normal University,
Beijing 100875, P. R. China}

\begin{abstract}
Using the framework first presented by Ruf and Sani in \cite{RufSani}, we give a proof of an Adams type inequality which can be applied to the functional
$$J_{\epsilon}(u) =\f{1}{2}\int_{\mathbb{R}^{2m}}\le(|\nabla^m u|^2+\sum_{\gamma=0}^{m-1}a_\gamma (x)|\nabla^\gamma u|^2\ri)dx
-\int_{\mathbb{R}^{2m}}\f{F(x,u)}{|x|^\beta}dx
-\epsilon\int_{\mathbb{R}^{2m}}hudx.$$
Under two kinds of assumptions on the nonlinearity, we estimate the min-max level of the functional.
As an application, a multiplicity result for the related singular quasilinear elliptic equation is proved.
\end{abstract}


\end{frontmatter}

\section{Introduction and main results}
Let $\nabla^\gamma u$, $\gamma\in \{0,1,2,\cdots,m\}$, be the $\gamma$-th order gradient of a function $u\in W^{m,2}(\mathbb{R}^{2m})$ which is defined by
$$\nabla^\gamma u:=\le\{\begin{array}{ll}
   \Delta^{\f{\gamma}{2}}u & \gamma\ \  \text{even},\\[1.5ex]
   \nabla\Delta^{\f{\gamma-1}{2}}u & \gamma\ \  \text{odd}.
   \end{array}\ri.$$
Here and throughout this paper, we use the notations that
$$\Delta^0 u=\nabla^0 u=u.$$
Consider the following nonlinear functional
\be\label{functional}
J_{\epsilon}(u) =\f{1}{2}\int_{\mathbb{R}^{2m}}\le(|\nabla^m u|^2+\sum_{\gamma=0}^{m-1}a_\gamma (x)|\nabla^\gamma u|^2\ri)dx
-\int_{\mathbb{R}^{2m}}\f{F(x,u)}{|x|^\beta}dx
-\epsilon\int_{\mathbb{R}^{2m}}hudx
\ee
which is related to the higher order partial differential equation
\be\label{equation}
(-\Delta)^m u+\sum_{\gamma=0}^{m-1}(-1)^\gamma\nabla^{\gamma}
\cdot(a_\gamma(x)\nabla^\gamma u)
=\f{f(x,u)}{|x|^\beta}+\epsilon h(x).
\ee
Here $m\geq 2$ is an even integer, $\ep$ is a small constant, the equation is defined on the whole Euclidean space of dimension $2m$, $0\leq\beta<2m$, $h(x)\nequiv 0$ belongs to the dual space of $E$ which will be defined later, $f(x,s):\mathbb{R}^{2m}\times\mathbb{R}\ra\mathbb{R}$ is a continuous function which satisfies some growth conditions and $a_\gamma(x)$ are continuous functions satisfying\vspace{.2cm}

\noindent${\bf (A_1)}$there exist positive constants $a_\gamma$, $\gamma=0,1,2,\cdots, m-1$, such that
$a_\gamma(x)\geq a_\gamma$ for all $x\in \mathbb{R}^{2m}$;

\noindent${\bf (A_2)}$ $(a_0(x))^{-1}\in L^1(\mathbb{R}^{2m})$.\vspace{.2cm}

This kind of equations has been extensively studied by many authors. When $m=1$, for the case $\beta=0$, the equation on a bounded domain $\Omega$ has been investigated in \cite{AdiYad, FDR, FMR, Zhao}. The corresponding $n$-Laplacian problem on a bounded domain also appears in many contexts, for example, in \cite{doo,Panda}. For an unbounded domain, the problem becomes different and for this case one can refer to \cite{Adimurthi, Cao, doo3} and the references therein. For the singular case, namely $0<\beta<n$, one can refer to \cite{AY, LL1, Yang1, Zhao1} and the references therein. Due to the the variational structure of this kind of equations, when $m=1$, usually the existence of solutions is related to the Moser-Trudinger type inequality. The inequality was first established by Truidinger \cite{Trudinger} and Moser \cite{Moser} and it says that, for a bounded domain $\Omega\subset \mathbb{R}^n$ and any $0\leq\alpha\leq \alpha_n=n\omega_{n-1}^{\f{1}{n-1}}$,
\be\label{trudinger}\sup_{u\in W_0^{1,n}(\Omega), \|\nabla u\|_{L^n(\Omega)}\leq 1}\int_{\Omega}e^{\alpha|u|^{\f{n}{n-1}}}dx<\infty,
\ee
where $\omega_{n-1}$ is the area of the unit sphere in $\mathbb{R}^n$.

When $m\geq 2$, related results about the corresponding higher order equations on bounded domains can be found in \cite{GazGruSqu, LazSch, LL3, ReiWet}. To deal with the higher order equations, we need a generalization of the Moser-Trudinger type inequality which is called the Adams type inequality. The classical Adams inequality given by Adams \cite{Adams} reads, for a bounded domain $\Omega\subset \mathbb{R}^n$ and any $0\leq\alpha\leq \alpha(m,n)$,
\be\label{adamsomega}
\sup_{u\in W_0^{m,\f{n}{m}}(\Omega), \|\nabla^m u\|_{L^{\f{n}{m}}(\Omega)}\leq 1}\int_{\Omega}e^{\alpha|u|^{\f{n}{n-m}}}dx<\infty,
\ee
where
$$\alpha(m,n)=\le\{\begin{array}{ll}
    \f{n}{\omega_{n-1}}\le(\f{\pi^{\f{n}{2}}2^m \Gamma\le(\f{m+1}{2}\ri)}{\Gamma\le(\f{n-m+1}{2}\ri)}\ri)^{\f{n}{n-m}}&\ \ m\ \ \text{odd},\\[1.5ex]
    \f{n}{\omega_{n-1}}\le(\f{\pi^{\f{n}{2}}2^m \Gamma\le(\f{m}{2}\ri)}{\Gamma\le(\f{n-m}{2}\ri)}\ri)^{\f{n}{n-m}}&\ \ m\ \ \text{even}.
   \end{array}\ri.$$
After Adams' work, many authors extended the inequality on a bounded domain from different points of view, for example, one can see \cite{Alb, Fon, Tarsi, YangZhao}. In particular, we mention the following singular Adams type inequality on a bounded domain \cite{LL2} which will be used later in our proof.\vspace{.2cm}

\noindent{\bf Theorem A} {\it Let $0\leq\beta<n$ and $\Omega\subset \mathbb{R}^n$ be a bounded domain. Then for any $0\leq\alpha\leq\le(1-\f{\beta}{n}\ri)\alpha(m,n)$, we have
\be\label{adamsbounded}
\sup_{u\in W_0^{m,\f{n}{m}}(\Omega),\|\nabla^m u\|_{L^{\f{n}{m}}(\Omega)}\leq 1}\int_{\Omega}\f{e^{\alpha |u|^{\f{n}{n-m}}}}{|x|^\beta}dx< \infty.
\ee
Moreover, when $m$ is an even number, the Sobolev space $W_0^{m,\f{n}{m}}(\Omega)$ in the above supremum can be replaced by the Sobolev space $W_{N}^{m,\f{n}{m}}(\Omega)$.}\vspace{.2cm}

In Theorem A, $W_{N}^{m,\f{n}{m}}(\Omega)$ is used to denote the space of functions with homogeneous Navier boundary conditions, namely,
$$W_{N}^{m,\f{n}{m}}(\Omega):=\le\{u\in W^{m,\f{n}{m}}(\Omega)\mid\Delta^\gamma u|_{\partial \Omega}=0\ \ \text{in the sense of traces for}\ \ 0\leq \gamma<\f{m}{2}\ri\}.$$
By definition, we have $W_0^{m,\f{n}{m}}(\Omega)\subset W_{N}^{m,\f{n}{m}}(\Omega)$,thus $W_{N}^{m,\f{n}{m}}(\Omega)$ is a larger Sobolev space.

It is easy to see that, when $\Omega\subseteq \mathbb{R}^n$ has infinite volume, the problem is that the integrals in both (\ref{trudinger}) and (\ref{adamsomega}) become infinite and the inequalities do not make sense. For the Moser-Trudinger type inequality, this problem was solved in \cite{Cao, Ruf} for dimension $n=2$ and in \cite{AdaTan, LR} for general dimension. Recently, for the Adams type inequality on an unbounded domain, Ruf and Sani \cite{RufSani} got the following result\vspace{.2cm}

\noindent{\bf Theorem B} {\it Let $m$ be an even integer
less than $n$ and
$\phi(t):=e^t-\sum_{\gamma=0}^{\gamma_{\f{n}{m}}-2}
\f{t^\gamma}{\gamma!}$, where
$\gamma_{\f{n}{m}}:=\min \le\{\gamma\in\mathbb{N}|\gamma\geq\f{n}{m}\ri\}\geq \f{n}{m}$.
There exists a constant $C_{m,n}>0$ such that, for any domain $\Omega\subseteq\mathbb{R}^n$,
\be\label{adamsRuf}
\sup_{u\in W_0^{m,\f{n}{m}}(\Omega),\|u\|_{m,n}\leq 1}\int_{\Omega}\phi\le(\alpha(m,n)|u|^{\f{n}{n-m}}\ri)dx\leq C_{m,n}
\ee
and this inequality is sharp.}\vspace{.2cm}

\noindent Hereafter we use $\|u\|_{m,n}$ to denote the norm of $u$ which is defined by
\bna
\|u\|_{m,n}:=\|(-\Delta+I)^{\f{m}{2}}u\|_{L^\f{n}{m}},
\ena
where $I$ denotes the identity operator.

After this, based on the ideas in Ruf and Sani's paper \cite{RufSani}, there are several generalizations of this result from different points of view. Lam and Lu \cite{LL4} improved Theorem B to the case that $m$ is an odd integer.
When $n=2m$ and $m\geq 2$ is an even integer, $\|u\|_{m,n}$ becomes
$$\|u\|_{m,2m}:=\|(-\Delta+I)^{\f{m}{2}}u\|_{L^2}.$$
But to be more suitable to use when considering equation (\ref{equation}), it is better to establish a singular Adams type inequality using the norm
$$\|u\|_{\tilde{E}}^2:=\int_{\mathbb{R}^{2m}}
\le(\sum_{\gamma=0}^m\tau_\gamma|\nabla^\gamma u|^2\ri)dx$$
instead of the norm $\|\cdot\|_{m,2m}$. Here $\tau_m=1$ and $\tau_\gamma>0$ for $\gamma=0,1,2,\cdots, m-1$.
For the nonsingular case, namely $\beta=0$, this was done in \cite{Yang} for $n=2m=4$ and in \cite{LL4} for general $n=2m$. When $0<\beta<n$, there are only results for the special dimension $n=2m=4$. In \cite{Yang}, Yang proved a result for the subcritical case $\alpha<\alpha(2,4)$ and in \cite{LL4}, Lam and Lu generalized the result to the critical case $\alpha=\alpha(2,4)$. In this paper, we consider the general case $n=2m$ and get the following theorem\vspace{.2cm}

\noindent{\bf Theorem 1.1} {\it Let $m\geq 2$ be an even integer, $\tau_m=1$, $\tau_\gamma>0$ for $\gamma=0,1,2,\cdots, m-1$ and $0\leq\beta<2m$, then for any $0\leq\alpha\leq\le(1-\f{\beta}{2m}\ri)\alpha(m,2m)$,
\be\label{adams}
\sup_{u\in W^{m,2}(\mathbb{R}^{2m}),\|u\|_{\tilde{E}}\leq 1}
\int_{\mathbb{R}^{2m}}
\f{e^{\alpha u^2}-1}{|x|^\beta}dx<\infty,
\ee
where $\alpha(m,2m)=(4\pi)^m m!$. Furthermore, the inequality is sharp.}\\

From now on we assume that $m\geq 2$ is an even integer and the dimension $n$ of the domain satisfies $n=2m$. Motivated by the Adams type inequality above, we assume the
following growth condition on the nonlinearity $f(x,s)$ of equation (\ref{equation}). \vspace{.2cm}

\noindent${\bf (H_1)}$ There exist constants $\alpha_0$, $b_1$, $b_2>0$ and $\theta\geq 1$ such that for all $(x,s)\in \mathbb{R}^{2m}\times\mathbb{R}$,
$$|f(x,s)|\leq b_1|s|+b_2|s|^\theta (e^{\alpha_0 s^2}-1).$$
\vspace{.2cm}

\noindent${\bf (H_2)}$ There exists $\mu>2$ such that for all $x\in\mathbb{R}^{2m}$ and $s\neq 0$,
$$0<\mu F(x,s)\equiv\mu\int_0^s f(x,t)dt\leq sf(x,s).$$
\vspace{.2cm}

\noindent${\bf (H_3)}$ There exist constants $R_0$, $M_0>0$ such that for all $x\in \mathbb{R}^{2m}$ and $|s|\geq R_0$,
$$0<F(x,x)\leq M_0 |f(x,s)|.$$
\vspace{.2cm}

Define a function space
$$E:=\le\{u\in W^{m,2}(\mathbb{R}^{2m}):
\int_{\mathbb{R}^{2m}}(|\nabla^m u|^2+\sum_{\gamma=0}^{m-1}a_\gamma(x)|\nabla^\gamma u|^2) dx<\infty\ri\}$$
and denote the norm of $u\in E$ by
$$\|u\|_E:=\le(\int_{\mathbb{R}^{2m}}
(|\nabla^m u|^2+\sum_{\gamma=0}^{m-1}a_\gamma(x)|\nabla^\gamma u|^2)dx\ri)^{\f{1}{2}}.$$
Here and in the sequel we use $E^*$ to denote the dual space of $E$ and assume $h(x)\in E^*$. Define a singular eigenvalue $\lambda_\beta$ by
\be\label{singulareigen}\lambda_\beta:=\inf_{u\in E\setminus \{0\}}\f{\|u\|_E^2}{\int_{\mathbb{R}^{2m}}\f{u^2}{|x|^\beta}dx}.
\ee
Moreover, we assume\vspace{.2cm}

\noindent${\bf (H_4)}$
$\limsup_{s\ra 0}\f{2|F(x,s)|}{s^2}<\lambda_\beta$ uniformly with respect to $x\in \mathbb{R}^{2m}$.
\vspace{.2cm}

The functional $J_\ep$ satisfies the geometric conditions of the mountain-pass theorem. The proof is similar to those in \cite{Yang} and \cite{Zhao1}. Namely, there exist two constant $r_\ep>0$ and $\vartheta_\ep>0$ such that $J_\ep(u)\geq \vartheta_\ep$ when $\|u\|_E=r_\ep$ and there exists some $e\in E$ satisfying $\|e\|_E>r_\ep$ such that  $J_\ep(e)<0$. Moreover, $J_\ep(0)=0$. Then the min-max level $C_M$ of $J_\ep$ is defined by
$$C_M=\min_{l\in \mathcal{L}}\max_{u\in l}J_\ep(u),$$
where $\mathcal{L}=\{l\in \mathcal{C}([0,1],E):l(0)=0,l(1)=e\}$. It is obvious that $C_M$ has a lower bound $\vartheta_\ep$, namely
$C_M \geq \vartheta_\ep$. We also want to give an explicit upper bound of $C_M$. To this end, we need the following additional assumptions\vspace{.2cm}

\noindent${\bf (H_5)}$ $\liminf_{s\ra +\infty}sf(x,s)e^{-\alpha_0 s^2}=+\infty$ uniformly with respect to $x\in \mathbb{R}^{2m}$
\vspace{.2cm}

\noindent or\vspace{.2cm}

\noindent${\bf (H_5)'}$ There exist constants $p>2$ and $C_p$ such that
$$|f(s)|\geq C_p |s|^{p-1},$$
where
$$C_p>\le(\f{p-2}{p}\ri)^{\f{p-2}{2}}
\le(\f{\alpha_0}{\le(1-\f{\beta}{2m}\ri)(4\pi)^m m!}\ri)^{\f{p-2}{2}}S_p^p,$$
\be\label{sp}
S_p^p:=\inf_{u\in E\setminus \{0\}}\f{\|u\|_E}
{\le(\int_{\mathbb{R}^{2m}}\f{u^p}{|x|^\beta}dx\ri)^{\f{1}{p}}}.
\ee
\vspace{.2cm}

Under each of these two assumptions, we can get the same estimate on the min-max level of (\ref{functional}).
Precisely, we have\vspace{.2cm}

\noindent{\bf Theorem 1.2} {\it Assume either $(H_5)$ or $(H_5)'$, together with $(H_2)$ and $(H_3)$, then there exists $\ep_0>0$ such that, for any $0<\ep\leq \ep_0$, the min-max level $C_M$ of (\ref{functional}) satisfies
\be\label{minmaxlevel}C_M< \le(1-\f{\beta}{2m}\ri)\f{(4\pi)^m m!}{2\alpha_0}.
\ee}\\

We remark that the above two assumptions on $f(x,s)$ can not cover each other. For details, one can refer to \cite{Yang2, Zhao1} for examples of $f(x,s)$ which can not satisfy these two assumptions simultaneously. As an application of the above estimate, we can get the following multiplicity result of equation (\ref{equation}). We can see later that the estimate on $C_M$ plays a crucial role in the proof of Theorem 1.3.\vspace{.2cm}

\noindent{\bf Theorem 1.3} {\it Assume either $(H_5)$ or $(H_5)'$, together with $(H_1)-(H_4)$, then there exists $\ep_1>0$ such that, for any $0<\ep\leq \ep_1$, the equation (\ref{equation}) has at least two distinct weak solutions.}\\

We organize this paper as follows: In Section 2, we prove the Adams type inequality (Theorem 1.1). In Section 3, we estimate the min-max level of functional (\ref{functional}) (Theorem 1.2). As an application of these two theorems, we prove the multiplicity result in Section 4 (Theorem 1.3).

\section{Adams type inequality}

Before the proof of Theorem 1.1, we point out that for $\tau_\gamma\leq a_\gamma$, $\gamma=0,1,2,\cdots, m-1$,
$$\|u\|_{\tilde{E}}\leq \|u\|_E.$$
This is the reason why Theorem 1.1 can be used in the study of equation (\ref{equation}). Using the Sobolev norm
$$\|u\|_{W^{m,2}}:=\le(\sum_{\gamma=0}^m
\|\nabla^\gamma u\|_{L^2}^2\ri)^{1/2},$$
it is easy to see that the norm $\|\cdot\|_{\tilde{E}}$ is equivalent to the norm $\|\cdot\|_{W^{m,2}}$. Another fact worth to emphasize is the following lemma\vspace{.2cm}

\noindent{\bf Lemma 2.1.} {\it Under assumptions $(A_1)$ and $(A_2)$, we have that the space $E$ is compactly embedded into the space $L^q(\mathbb{R}^{2m})$ for any $q\geq 1$.}\\

The proof of this lemma is essentially the same as the proof of Lemma 3.6 in \cite{Yang}. But for the convenience of readers, we give a proof here.\vspace{.2cm}

\noindent{\bf Proof.} When $q\geq 2$, it is easy to see that the embedding $E\hookrightarrow L^q(\mathbb{R}^{2m})$ is continuous. When $q=1$, H\"{o}lder's inequality and $(A_2)$ imply that
$$\int_{\mathbb{R}^{2m}}|u|dx\leq \le(\int_{\mathbb{R}^{2m}}\f{1}{a_0(x)}dx\ri)^{\f{1}{2}}
\le(\int_{\mathbb{R}^{2m}}a_0(x)u^2dx\ri)^{\f{1}{2}}
\leq \le(\int_{\mathbb{R}^{2m}}\f{1}{a_0(x)}dx\ri)^{\f{1}{2}}\|u\|_E.$$
When $1<q<2$, we have
$$\int_{\mathbb{R}^{2m}}|u|^q dx
\leq \int_{\mathbb{R}^{2m}}(|u|+u^2) dx
\leq \le(\int_{\mathbb{R}^{2m}}\f{1}{a_0(x)}dx\ri)^{\f{1}{2}}\|u\|_E
+\f{1}{a_0}\|u\|_E^2.$$
Thus we have that, for any $q\geq 1$, the embedding $E\hookrightarrow L^q(\mathbb{R}^{2m})$ is continuous.

Next we prove that the embedding is also compact. Suppose $\{u_k\}\subset E$ is a bounded sequence, we need to prove that $u_k$ converges to some $u\in E$ strongly in $L^q(\mathbb{R}^{2m})$ up to a subsequence for any $q\geq 1$.

$(A_2)$ implies that, for any $\epsilon>0$, there exists $R_\epsilon>0$ such that
$$\int_{|x|>R_\ep} \f{1}{a_0(x)} dx<\epsilon^2.$$
Since $\{u_k\}$ is a bounded sequence, up to subsequence, we can assume that $u_k$ converges to some $u$ strongly in $L^1(B_{R_\ep})$. When $q=1$, we have
\bea\label{q1}
\int_{\mathbb{R}^{2m}}|u_k-u| dx
&=&\int_{|x|\leq R_\ep}|u_k-u| dx
+\int_{|x|>R_\ep}|u_k-u| dx\nonumber\\
&\leq& \int_{|x|\leq R_\ep}|u_k-u| dx
+\le(\int_{|x|>R_\ep}\f{1}{a_0(x)}dx\ri)^{\f{1}{2}}
\le(\int_{|x|>R_\ep}a_0(x)|u_k-u|^2dx\ri)^{\f{1}{2}}\nonumber\\
&\leq&\int_{|x|\leq R_\ep}|u_k-u| dx
+\epsilon \|u_k-u\|_E.
\eea
Noticing that $\ep$ can be arbitrarily small, we get from (\ref{q1}) that
$$\lim_{k\ra \infty}\int_{\mathbb{R}^{2m}}|u_k-u|dx=0.$$

When $q>1$, we have
\bea\label{q2}
\int_{\mathbb{R}^{2m}}|u_k-u|^q dx
&\leq&\le(\int_{\mathbb{R}^{2m}}|u_k-u|dx\ri)^{\f{1}{2}}
\le(\int_{\mathbb{R}^{2m}}|u_k-u|^{2q-1}dx\ri)^{\f{1}{2}}\nonumber\\
&\leq&C\le(\int_{\mathbb{R}^{2m}}|u_k-u|dx\ri)^{\f{1}{2}}\ra 0
\eea
as $k\ra \infty$. Here we used the continuous embedding $E\hookrightarrow L^{2q-1}(\mathbb{R}^{2m})$.
$\hfill\Box$\\

We remark here that the singular eigenvalue $\lambda_\beta$ defined in (\ref{singulareigen}) and $S_p$ defined in (\ref{sp}) are both positive constants for any $0\leq \beta<2m$. When $\beta=0$, $(A_1)$ gives us that $\lambda_0\geq a_0>0$. When $0<\beta<2m$, we have
\bna
\int_{\mathbb{R}^{2m}}\f{u^2}{|x|^\beta}dx&\leq& \int_{|x|>1}u^2dx+\le(\int_{|x|\leq 1}|u|^{2q}dx\ri)^{\f{1}{q}}
\le(\int_{|x|\leq 1}\f{1}{|x|^{\beta q'}}dx\ri)^{\f{1}{q'}}\\
&\leq&C\|u\|_E^2,
\ena
where $\f{1}{q}+\f{1}{q'}=1$ and $0<\beta q'<2m$. This implies that $\lambda_\beta\geq \f{1}{C}>0$. Similarly, we can prove $S_p>0$.

To prove Theorem 1.1, we first give several definitions. Let $B_R$ be an open ball centered at $0$ with radius $R>0$. If $u: B_R\ra \mathbb{R}$ is a measurable function, the distribution function of $u$ is defined by
$$\mu_u(t):=\mathcal{M}\le(\{x\in B_R\ |\ |u(x)|>t\}\ri)\ \ \ \forall t\geq 0,$$
where $\mathcal{M}(\cdot)$ denotes the Lebesgue measure of a set in $\mathbb{R}^n$. The decreasing rearrangement of $u$ is defined by
$$u^\star(s):=\inf\{t\geq 0\ |\ \mu_u(t)<s\}\ \ \forall s\in [0, \mathcal{M}(B_R)].$$
Finally, the spherically symmetric decreasing rearrangement of $u$ is defined by
$$u^*(x):=u^\star(\sigma_n|x|^n)\ \ \forall x\in B_R,$$
where $\sigma_n$ is the volume of the unit ball in $\mathbb{R}^n$.

Now we begin to prove Theorem 1.1 by using the framework of Ruf and Sani's work \cite{RufSani}. After \cite{RufSani}, similar ideas were also used in \cite{LL4}, \cite{LL2} and \cite{Yang} .\vspace{.2cm}

\noindent{\bf Proof of Theorem 1.1.} For any $u\in W^{m,2}(\mathbb{R}^{2m})$ and $\tilde{\rho}>0$, direct computations give that
\be\label{expansion}
\int_{\mathbb{R}^{2m}} \le((-\Delta+\tilde{\rho} I)^\f{m}{2} u\ri)^2dx=\sum_{\gamma=0}^m C(m,\gamma) \tilde{\rho}^{m-\gamma}\int_{\mathbb{R}^{2m}}|\nabla^\gamma u|^2dx,
\ee
where $C(m,\gamma)=\f{m!}{\gamma!(m-\gamma)!}$. In particular, $C(m,m)=C(m,0)=1$ and $C(m,1)=m$. Define $\rho_\gamma=\le(\f{\tau_\gamma}{C(m,\gamma)}\ri)^{\f{1}{m-\gamma}}$, $\gamma=\{0,1,2,\cdots,m\}$ and let $\rho=\min\{\rho_0,\rho_1,\cdots,\rho_m\}$. Then (\ref{expansion}) tells us that
$$\int_{\mathbb{R}^{2m}} \le((-\Delta+\rho I)^\f{m}{2} u\ri)^2dx\leq \|u\|_{\tilde{E}}^2.$$
So if we can prove that
$$\sup_{u\in W^{m,2}(\mathbb{R}^{2m}),\int_{\mathbb{R}^{2m}} \le((-\Delta+\rho I)^\f{m}{2} u\ri)^2dx\leq 1}
\int_{\mathbb{R}^{2m}}
\f{e^{\alpha u^2}-1}{|x|^\beta}dx<\infty,$$
the inequality in Theorem 1.1 is proved immediately.

By density of $C_0^\infty (\mathbb{R}^{2m})$ in $W^{m,2}(\mathbb{R}^{2m})$, we can find a sequence $\{u_k\}\subset C_0^\infty (\mathbb{R}^{2m})$ such that $u_k\ra u$ in $W^{m,2}(\mathbb{R}^{2m})$. Without loss of generality, we can assume that $\int_{\mathbb{R}^{2m}} \le((-\Delta+\rho I)^\f{m}{2} u_k\ri)^2dx=1$, for otherwise we can use $\tilde{u}_k=\f{u_k}{\le(\int_{\mathbb{R}^{2m}} \le((-\Delta+\rho I)^\f{m}{2} u_k\ri)^2dx\ri)^{\f{1}{2}}}$ instead of $u_k$.

Suppose, for a fixed $k$, supp $u_k\subset B_{R_k}$. Define $$f_k=(-\Delta+\rho I)^\f{m}{2} u_k$$
and use $f_k^*$ to denote the spherically symmetric decreasing rearrangement of $f_k$. Consider the equation
$$\le\{\begin{array}{ll}
   (-\Delta+\rho I)^\f{m}{2} v_k=f_k^* &\ \ \text{in}\ \  B_{R_k},\\[1.5ex]
    v_k\in W_{N}^{m,2}(B_{R_k}).
   \end{array}\ri.$$
By properties of rearrangement (see \cite{Kav, LL2, RufSani}), we have
\be\label{muk1}\int_{B_{R_k}}\le((-\Delta+\rho I)^\f{m}{2} v_k\ri)^2dx=\|f_k^*\|_{L^2(B_{R_k})}^2=\|f_k\|_{L^2(B_{R_k})}^2=
\int_{B_{R_k}}\le((-\Delta+\rho I)^\f{m}{2} u_k\ri)^2dx=1
\ee
and
\be\label{ukmuk}
\int_{B_{R_k}}\f{e^{\alpha u_k^2}-1}{|x|^\beta}dx\leq \int_{B_{R_k}}\f{e^{\alpha v_k^2}-1}{|x|^\beta}dx.
\ee

Let $r_0\geq 1$ be a constant to be determined later. If $R_k\leq r_0$, since
$$\|\nabla^m v_k\|_{L^2(B_{R_k})}^2\leq \int_{B_{R_k}}\le((-\Delta+\rho I)^\f{m}{2} v_k\ri)^2dx=1,$$
by Theorem A, we get
$$\int_{B_{R_k}}\f{e^{\alpha v_k^2}-1}{|x|^\beta}dx<C_{m,r_0},$$
where $C_{m,r_0}$ is some constant depending on $m$ and $r_0$ but not depending on $k$.

If $R_k>r_0$, rewrite $\int_{B_{R_k}}\f{e^{\alpha v_k^2}-1}{|x|^\beta}dx$ into
$$\int_{B_{r_0}}\f{e^{\alpha v_k^2}-1}{|x|^\beta}dx+\int_{B_{R_k}\setminus B_{r_0}}\f{e^{\alpha v_k^2}-1}{|x|^\beta}dx:=I_1+I_2.$$

Firstly, we estimate $I_1$. Define, for $\gamma=\{1,2,\cdots,\f{m}{2}\}$ and $x\in B_{r_0}$,
$$\xi_\gamma(|x|)=|x|^{m-2\gamma}.$$
Let $$g_k(x)=\sum_{\gamma=1}^{\f{m}{2}}d_{k,\gamma}\xi_\gamma(|x|),$$
where
$$d_{k,\gamma}=\f{\Delta^{\f{m}{2}-\gamma}v_k(r_0)
-\sum_{\eta=1}^{\gamma-1}d_{k,\eta}\Delta^{\f{m}{2}-\gamma}\xi_\eta(r_0)}
{\Delta^{\f{m}{2}-\gamma}\xi_\gamma(r_0)}.$$
Denote $\le(v_k(x)-g_k(x)\ri)$ by $\mu_k(x)$. By construction, we have (see Lemma 4.3 in \cite{RufSani}) $\nabla^m\mu_k=\nabla^m\nu_k$ in $B_{r_0}$, $\mu_k$ is a radial function with homogeneous Navier boundary conditions and for $0<|x|\leq r_0$,
\be\label{radialestimate}
v_k^2(x)\leq
\mu_k^2(x)\le(1+C_m\sum_{\gamma=1}^{\f{m}{2}}
r_0^{1-4\gamma}\|\Delta^{\f{m}{2}-\gamma}v_k\|_{W^{1,2}(B_{r_0})}^2
\ri)^2+C_{m,r_0}.
\ee
Since $r_0>1$, (\ref{radialestimate}) implies that 
\be\label{radialestimate1}
v_k^2(x)\leq
\mu_k^2(x)\le(1+C_m r_0^{-3}\sum_{\gamma=0}^{m-1}
\|\nabla^{\gamma}v_k\|_{L^2(B_{r_0})}^2
\ri)^2+C_{m,r_0}.
\ee
Define
$$\tilde{\mu}_k(x):=\mu_k(x)\le(1+C_m r_0^{-3}\sum_{\gamma=0}^{m-1}
\|\nabla^{\gamma}v_k\|_{L^2(B_{r_0})}^2
\ri)$$
and
$$C_{\min}:=\min\{C(m,\gamma)\rho^{m-\gamma}\ |\ 0\leq\gamma\leq m-1 \}.$$
Then we have
\bea\label{mutilde}
\|\nabla^m\tilde{\mu}_k\|_{L^2(B_{r_0})}^2&=&\le(1+C_m r_0^{-3}\sum_{\gamma=0}^{m-1}
\|\nabla^{\gamma}v_k\|_{L^2(B_{r_0})}^2\ri)^2
\|\nabla^m \mu_k\|_{L^2(B_{r_0})}^2\nonumber\\
&=&\le(1+C_m r_0^{-3}\sum_{\gamma=0}^{m-1}
\|\nabla^{\gamma}v_k\|_{L^2(B_{r_0})}^2\ri)^2
\|\nabla^m \nu_k\|_{L^2(B_{r_0})}^2\nonumber\\
&=&\le(1+C_m r_0^{-3}\sum_{\gamma=0}^{m-1}
\|\nabla^{\gamma}v_k\|_{L^2(B_{r_0})}^2\ri)^2
\le(1-\sum_{\gamma=0}^{m-1}C(m,\gamma)\rho^{m-\gamma}\|\nabla^\gamma v_k\|_{L^2{(B_{r_0})}}^2\ri)\nonumber\\
&\leq&\le(1+C_m r_0^{-3}\sum_{\gamma=0}^{m-1}
\|\nabla^{\gamma}v_k\|_{L^2(B_{r_0})}^2\ri)^2
\le(1-C_{\min}\sum_{\gamma=0}^{m-1}\|\nabla^\gamma v_k\|_{L^2{(B_{r_0})}}^2\ri),
\eea
where we have used (\ref{expansion}) and (\ref{muk1}) at the third equality. Choose $r_0^3\geq\max\{1, \f{C_m}{C_{\min}}\}$, we get from (\ref{mutilde}) that
\be\label{mutilde1}
\|\nabla^m\tilde{\mu}_k\|_{L^2(B_{r_0})}^2\leq 1.
\ee
Now by (\ref{radialestimate1}), we have
\bna
I_1&\leq&\int_{B_{r_0}}\f{e^{\alpha\mu_k^2(x)\le(1+C_m r_0^{-3}\sum_{\gamma=0}^{m-1}
\|\nabla^{\gamma}v_k\|_{L^2(B_{r_0})}^2
\ri)^2+\alpha C_{m,r_0}}-1}{|x|^\beta} dx\\
&=&\int_{B_{r_0}}\f{e^{\alpha\tilde{\mu}_k^2(x)
+\alpha C_{m,r_0}}-1}{|x|^\beta} dx\\
&\leq&e^{\alpha C_{m,r_0}}\int_{B_{r_0}}\f{e^{\alpha\tilde{\mu}_k^2(x)
}}{|x|^\beta} dx.
\ena
Then (\ref{mutilde1}) and Theorem A imply that 
\be\label{i1}
I_1\leq C_{m,r_0,\alpha,\rho}.
\ee

Secondly, we deal with $I_2$. The radial lemma (see Chapter 6 in \cite{Kav}) tells us that for $v_k\in W_N^{m,2}(B_{R_k})\subset W^{1,2}(\mathbb{R}^{2m})$, we have
\be\label{radiallemma}
|v_k(x)|\leq\sqrt{\f{(m-1)!}{\pi^m}}{|x|^{\f{1}{2}-m}}
\|v_k\|_{W^{1,2}(\mathbb{R}^{2m})}\ \ \text{a.e. in}\ \ \mathbb{R}^{2m}.
\ee
Take $r_0^{2m-1}\geq \f{(m-1)!}{\pi^m}
\f{1}{\min\{m\rho^{m-1},\rho^{m}\}}$. If $|x|\geq r_0$, by (\ref{radiallemma}), we have
\be\label{v1}
|v_k(x)|\leq \sqrt{\min \{m\rho^{m-1},\rho^{m}\}}\|v_k\|_{W^{1,2}(\mathbb{R}^{2m})}.
\ee
On the other hand, we have
\bea\label{v2}
\min \{m\rho^{m-1},\rho^{m}\}\|v_k\|^2
_{W^{1,2}(\mathbb{R}^{2m})}&\leq&\int_{B_{R_k}}
\le(m\rho^{m-1}|\nabla v_k|^2+\rho^{m}v_k^2\ri)dx\nonumber\\
&\leq&\int_{B_{R_k}}
\le((-\Delta+\rho I)^\f{m}{2} v_k\ri)^2dx\nonumber\\
&=&1.
\eea
Obviously, (\ref{v1}) and (\ref{v2}) imply that $|v_k(x)|\leq 1$ for any $|x|\geq r_0$. It follows that
\bea\label{i2}
I_2&\leq&\f{1}{r_0^\beta}\int_{B_{R_k}\setminus B_{r_0}}(e^{\alpha v_k^2}-1)dx\nonumber\\
&=&\f{1}{r_0^\beta}\int_{B_{R_k}\setminus B_{r_0}}\sum_{l=1}^\infty \f{\alpha^l v_k^{2l}}{l!}dx\nonumber\\
&\leq&\f{1}{r_0^\beta}\int_{B_{R_k}\setminus B_{r_0}}\sum_{l=1}^\infty \f{\alpha^l v_k^{2}}{l!}dx\nonumber\\
&\leq&\f{1}{r_0^\beta}\sum_{l=1}^\infty \f{\alpha^l }{l!}\|v_k\|_{W^{1,2}(\mathbb{R}^{2m})}^2\nonumber\\
&\leq&\f{1}{\min \{m\rho^{m-1},\rho^{m}\}r_0^\beta}\sum_{l=1}^\infty \f{\alpha^l }{l!}\nonumber\\
&\leq&C_{m,r_0,\alpha,\rho}.
\eea
Take $r_0\geq \max\{1, \le(\f{C_m}{C_{\min}}\ri)^{\f{1}{3}}, \le(\f{(m-1)!}{\pi^m}
\f{1}{\min\{m\rho^{m-1},\rho^{m}\}}\ri)^{\f{1}{2m-1}}\}$. Then Fatou's lemma together with (\ref{ukmuk}), (\ref{i1}) and (\ref{i2}) tells that
$$\int_{\mathbb{R}^{2m}}
\f{e^{\alpha u^2}-1}{|x|^\beta}dx\leq \liminf_{k\ra \infty}\int_{\mathbb{R}^{2m}}
\f{e^{\alpha u_k^2}-1}{|x|^\beta}dx\leq C_{m,r_0,\alpha,\rho}$$
and the proof of the inequality is finished.

To prove the sharpness of the inequality, we need a sequence of test functions. For this reason,  we postpone the proof of sharpness till the end of Section 3.
$\hfill\Box$\\

\section{Min-max level}

In this section, we estimate the min-max level of $J_\epsilon$. Firstly, we define a sequence of functions $\tilde{\phi}_k(x)$ by
$$\tilde{\phi}_k(x)=\le\{\begin{array}{lll}
   \sqrt{\f{\log k}{2M}}+\f{1}{\sqrt{2M\log k}}\sum_{\gamma=1}^{m-1}\f{(1-k |x|^2)^\gamma}{\gamma} & |x|\in [0,\f{1}{\sqrt{k}}),\\[1.5ex]
   -\sqrt{\f{2}{M\log k}}\log |x| & |x|\in [\f{1}{\sqrt{k}},1),\\[1.5ex]
   \zeta_k(x) & |x|\in [1,\infty)
   \end{array}\ri.$$
where
$$M=\f{(4\pi)^m (m-1)!}{2}, \zeta_k\in C_0^\infty(B_2(0)), \zeta_k\mid_{\partial B_1(0)}=\zeta_k\mid_{\partial B_2(0)}=0.$$
Moreover, for $\gamma=\{1,2,\cdots, m-1\}$, $\f{\partial^\gamma\zeta_k}{\partial r^\gamma}\mid_{\partial B_1(0)}=(-1)^\gamma(\gamma-1)!\sqrt{\f{2}{M\log k}}$, $\f{\partial^\gamma\zeta_k}{\partial r^\gamma}\mid_{\partial B_2(0)}=0$  and $\zeta_k$, $|\nabla^\gamma\zeta_k|$ and $|\nabla^m\zeta_k|$ are all $O\le(\f{1}{\sqrt{\log k}}\ri)$. Obviously, $\tilde{\phi}_k(x)$ are continuous functions defined on $\mathbb{R}^{2m}$ with compact supports .

To estimate the $W^{m,2}$ norms of $\tilde{\phi}_k$, we need the following lemma.\vspace{.2cm}

\noindent{\bf Lemma 3.1.} {\it For $\gamma=\{1,2,\cdots, m-1\}$, the $\gamma$-th order derivatives of $\tilde{\phi}_k$ with respect to $r$ satisfy
\be\label{lr1}
\lim_{r\ra {\f{1}{\sqrt{k}}}^-}\f{\partial^\gamma\tilde{\phi}_k}{\partial r^\gamma}=\lim_{r\ra {\f{1}{\sqrt{k}}}^+}\f{\partial^\gamma\tilde{\phi}_k}{\partial r^\gamma}
\ee
and
\be\label{lr2}
\lim_{r\ra {1}^-}\f{\partial^\gamma\tilde{\phi}_k}{\partial r^\gamma}=\lim_{r\ra {1}^+}\f{\partial^\gamma\tilde{\phi}_k}{\partial r^\gamma},
\ee
where $r=|x|$.}\\

\noindent{\bf Proof.} Direct computations give that, when
$\f{1}{\sqrt{k}}\leq r<1$, the $\gamma$-th order derivatives of
$\tilde{\phi}_k$ with respect to $r$ are
$$(-1)^\gamma(\gamma-1)!\sqrt{\f{2}{M\log k}}r^{-\gamma}.$$
Combining this with our assumptions on $\zeta_k$, we get (\ref{lr2}).

To get (\ref{lr1}), we consider the following functions of $r$
$$t_k(r)=-\f{1}{\sqrt{2M\log k}}\log r.$$
The Taylor series of $t_k(r)$ at $\f{1}{k}$ is
\bna
t_k(r)&=&\sum_{\gamma=0}^{\infty}\f{t_k^{(\gamma)}(\f{1}{k})}
{\gamma !}(r-\f{1}{k})^\gamma\\
&=&\sqrt{\f{\log k}{2M}}+\f{1}{\sqrt{2M\log k}}\sum_{\gamma=1}^{\infty}\f{(1-k r)^\gamma}{\gamma}.
\ena
We use $\tilde{t}_k(r)$ to denote the summation of the first $m$ terms of the series, namely,
$$\tilde{t}_k(r)=\sqrt{\f{\log k}{2M}}+\f{1}{\sqrt{2M\log k}}\sum_{\gamma=1}^{m-1}\f{(1-k r)^\gamma}{\gamma}.$$
It is easy to know that, at $r=\f{1}{k}$, for $\gamma=\{1,2,\cdots,m-1\}$, the $\gamma$-th order derivatives of $t_k(r)$ equal to those of $\tilde{t}_k(r)$ respectively. By the definitions of $\tilde{\phi}_k$, we have
$$\tilde{\phi}_k(x)=\le\{\begin{array}{ll}
    \tilde{t}_k(r^2) & r\in [0,\f{1}{\sqrt{k}}),\\[1.5ex]
    t_k(r^2)& r\in [\f{1}{\sqrt{k}},1).
  \end{array}\ri.$$
This fact implies (\ref{lr1}) immediately.
$\hfill\Box$\\

We remark that to find the extremal of Adams inequality, Adams has constructed a sequence of functions in \cite{Adams} which has properties similar to our sequence. But at first, Adams' functions have no explicit expressions. Moreover, our functions are defined on the whole space $\mathbb{R}^{2m}$ instead of a bounded domain $\Omega\subset \mathbb{R}^{2m}$.

We claim that $\tilde{\phi}_k(x)\in W_0^{m,2}(\mathbb{R}^{2m})$ and, for $\gamma=\{0,1,\cdots,m-1\}$,
$$\|\nabla^\gamma\tilde{\phi}_k\|_{L^2}^2=O\le(\f{1}{\log k}\ri),$$
while
$$\|\nabla^m\tilde{\phi}_k\|_{L^2}^2
=\|\Delta^{\f{m}{2}}\tilde{\phi}_k\|_{L^2}^2=1+O\le(\f{1}{\log k}\ri).$$

To prove the claim, we first point out that, by Lemma 3.1 and the formula for integration by parts, we can get the weak derivatives of $\tilde{\phi}_k(x)$ until order $m$ by computations on each part of the domain. Therefore, we can estimate the $W^{m,2}$ norms of $\tilde{\phi}_k$ respectively.\vspace{.2cm}

\noindent{\bf Part I $\mathbb{R}^{2m}\setminus B_1(0)$.}

By definitions, since $\zeta_k\in C_0^\infty(B_2(0))$, we have, on $\mathbb{R}^{2m}\setminus B_1(0)$,
\be\label{norm11}
\|\tilde{\phi}_k\|_{L^2}^2=\|\zeta_k\|_{L^2}^2=
O\le(\f{1}{\log k}\ri)
\ee
and, for $\gamma=\{1,2,\cdots, m\}$,
\be\label{norm12}
\|\nabla^\gamma\tilde{\phi}_k\|_{L^2}^2=
\|\nabla^\gamma\zeta_k\|_{L^2}^2=
O\le(\f{1}{\log k}\ri).
\ee
\vspace{.2cm}

\noindent{\bf Part II $B_1(0)\setminus B_{\f{1}{\sqrt{k}}}(0)$.}

When $\gamma=1$, it is easy to get
\be\label{gamma1}|\nabla^\gamma\tilde{\phi}_k|=|\nabla \tilde{\phi}_k|=-\sqrt{\f{2}{M\log k}}r^{-1}.
\ee
For higher order derivatives, noticing the fact that, for any integer $1\leq l\leq \f{m}{2}$,
$$\Delta^l \log r=(-1)^{l-1} 2^{2l-1} \f{(l-1)!(m-1)!}{(m-l-1)!}r^{-2l},$$
we have, when $\gamma$ is odd and $1<\gamma<m$,
\bea\label{gamma1m}
|\nabla^\gamma\tilde{\phi}_k|&=&\le|\f{\partial}{\partial r}(\Delta^{\f{\gamma-1}{2}}\tilde{\phi}_k)\ri|\nonumber\\
&=&\le|(-1)^{\f{\gamma+1}{2}}
2^{\gamma-2}(\gamma-1)\sqrt{\f{2}{M\log k}}\f{(m-1)!\le(\f{\gamma-3}{2}\ri)!}{\le(m-\f{\gamma+1}{2}\ri)!}
r^{-(\gamma+1)}
\vec{x}\ri|\nonumber\\
&=&\le|2^{\gamma-2}(\gamma-1)\sqrt{\f{2}{M\log k}}\f{(m-1)!\le(\f{\gamma-3}{2}\ri)!}{\le(m-\f{\gamma+1}{2}\ri)!}r^{-\gamma}\ri|.
\eea
When $\gamma$ is even and $2\leq\gamma\leq m$,
\be\label{gamma2m}\nabla^\gamma\tilde{\phi}_k=
\Delta^{\f{\gamma}{2}}\tilde{\phi}_k=
(-1)^{\f{\gamma}{2}}
2^{\gamma-1}\sqrt{\f{2}{M\log k}}\f{(m-1)!\le({\f{\gamma}{2}}-1\ri)!}{\le(m-{\f{\gamma}{2}}-1\ri)!}r^{-\gamma}.
\ee
In particular, we have
$$\nabla^m\tilde{\phi}_k=\Delta^{\f{m}{2}}\tilde{\phi}_k=
(-1)^{\f{m}{2}}
2^{m-1}\sqrt{\f{2}{M\log k}}(m-1)!r^{-m},$$
which gives us that, on $B_1(0)\setminus B_{\f{1}{\sqrt{k}}}(0)$,
\bna
\|\nabla^m\tilde{\phi}_k\|_{L^2}^2&=&\int_{B_1(0)\setminus B_{\f{1}{\sqrt{k}}}(0)} \f{2^{2m-1}\le((m-1)!\ri)^2}{M\log k}r^{-2m}dx\\
&=&\f{2^{2m-1}\le((m-1)!\ri)^2}{M\log k}\int_{\f{1}{\sqrt{k}}}^1\omega_{2m-1}r^{-1}dr\\
&=&\f{4^{m-1}\le((m-1)!\ri)^2\omega_{2m-1}}{M}.
\ena
Substituting $M=\f{(4\pi)^m (m-1)!}{2}$ and $\omega_{2m-1}=\f{2\pi^m}{(m-1)!}$, we get
\be\label{normm}
\|\nabla^m\tilde{\phi}_k\|_{L^2}^2=1.
\ee
When $\gamma=0$, by integrating by parts and the definitions of $\tilde{\phi}_k$, we get
\bea\label{norm0}
\|\tilde{\phi}_k\|_{L^2}^2&=&\int_{B_1(0)\setminus B_{\f{1}{\sqrt{k}}}(0)} \f{2\log^2 r}{M\log k}dx\nonumber\\
&=&\f{\omega_{2m-1}}{M\log k}\le(\f{1}{2m^3}-\f{\log^2 k}{4m k^{m}}-\f{\log k}{2m^2 k^{m}}-\f{1}{2m^3 k^{m}}\ri)\nonumber\\
&=&O\le(\f{1}{\log k}\ri).
\eea
Similarly, when $\gamma=1$, (\ref{gamma1})gives us that
\bea\label{norm1}
\|\nabla\tilde{\phi}_k\|_{L^2}^2
&=&\int_{B_1(0)\setminus B_{\f{1}{\sqrt{k}}}(0)} \f{2r^{-2}}{M\log k}dx\nonumber\\
&=&\f{\omega_{2m-1}}{(m-1)M\log k}\le(1-\f{1}{k^{m-1}}\ri)\nonumber\\
&=&O\le(\f{1}{\log k}\ri).
\eea
When $\gamma$ is odd and $1<\gamma<m$, by (\ref{gamma1m}), we have
\bea\label{norm1m}
\|\nabla^\gamma\tilde{\phi}_k\|_{L^2}^2&=&\int_{B_1(0)\setminus B_{\f{1}{\sqrt{k}}}(0)} \f{(\gamma-1)^2 2^{2\gamma-3}}{M\log k}\le(\f{(m-1)!\le(\f{\gamma-3}{2}\ri)!}{\le(m-\f{\gamma+1}{2}\ri)!}\ri)^2
r^{-2\gamma}dx\nonumber\\
&=&\f{(\gamma-1)^2 4^{\gamma-2}\omega_{2m-1}}{(m-\gamma)M\log k}\le(\f{(m-1)!\le(\f{\gamma-3}{2}\ri)!}{\le(m-\f{\gamma+1}{2}\ri)!}\ri)^2
\le(1-\f{1}{k^{m-\gamma}}\ri)\nonumber\\
&=&O\le(\f{1}{\log k}\ri).
\eea
When $\gamma$ is even and $2\leq \gamma\leq m-2$, by (\ref{gamma2m}), we have
\bea\label{norm2m}
\|\nabla^\gamma\tilde{\phi}_k\|_{L^2}^2&=&\int_{B_1(0)\setminus B_{\f{1}{\sqrt{k}}}(0)} \f{2^{2\gamma-1}}{M\log k}\le(\f{(m-1)!\le(\f{\gamma}{2}-1\ri)!}{\le(m-\f{\gamma}{2}-1\ri)!}\ri)^2 r^{-2\gamma}dx\nonumber\\
&=&\f{4^{\gamma-1}\omega_{2m-1}}{M\log k(m-\gamma)}\le(\f{(m-1)!\le(\f{\gamma}{2}-1\ri)!}{\le(m-\f{\gamma}{2}-1\ri)!}\ri)^2
\le(1-\f{1}{k^{m-\gamma}}\ri)\nonumber\\
&=&O\le(\f{1}{\log k}\ri).
\eea\vspace{.2cm}

\noindent{\bf Part III $B_{\f{1}{\sqrt{k}}}(0)$.}

We have
$$\|\tilde{\phi}_k\|_{W^{m,2}}\leq \le\|\sqrt{\f{\log k}{2M}}\ri\|_{W^{m,2}}+\sum_{\gamma=1}^{m-1}
\le\|\f{k^\gamma}{\gamma\sqrt{2M\log k}}(\f{1}{k}-|x|^2)^\gamma\ri\|_{W^{m,2}}.$$
Direct computations show that
\bea\label{norm31}
\le\|\sqrt{\f{\log k}{2M}}\ri\|_{W^{m,2}(B_{\f{1}{\sqrt{k}}})}^2
&=&\int_{B_{\f{1}{\sqrt{k}}}}\f{\log k}{2M}dx\nonumber\\
&=&\f{\omega_{2m-1}\log k}{4mM k^{m}}\nonumber\\
&=&o\le(\f{1}{\log k}\ri).
\eea
Furthermore, by integrating by parts, we have
\bea\label{norm32}
\le\|\f{k^\gamma}{\gamma\sqrt{2M\log k}}(\f{1}{k}-|x|^2)^\gamma\ri\|_{W^{m,2}}^2
&=&
\sum_{\eta=0}^m\le\|\f{k^\gamma}{\gamma\sqrt{2M\log k}}\nabla^\eta\le((\f{1}{k}-|x|^2)^\gamma\ri)\ri\|_{L^2(B_{\f{1}{\sqrt{k}}}(0))}^2\nonumber\\
&=&\le\{\begin{array}{ll}
   \sum_{\eta=0}^{2\gamma} O\le(\f{1}{k^{m-\eta}\log k}\ri)\ \ \gamma<\f{m}{2},\\[2.5ex]
   \sum_{\eta=0}^m O\le(\f{1}{k^{m-\eta}\log k}\ri)\ \ \gamma\geq \f{m}{2}
   \end{array}\ri.\nonumber\\
&=&O\le(\f{1}{\log k}\ri).
\eea

Combining (\ref{norm11}), (\ref{norm12}) and (\ref{normm}-\ref{norm32}), we prove the claim that $$\|\tilde{\phi}_k\|_{W^{m,2}(\mathbb{R}^{2m})}^2=1
+O\le(\f{1}{\log k}\ri),$$
or equivalently, 
$$\|\tilde{\phi}_k\|_{E}^2=1
+O\le(\f{1}{\log k}\ri),$$

Define $$\phi_k(x)=\f{\tilde{\phi}_k(x)}{\|\tilde{\phi}_k\|_{E}}.$$ We have $\|\phi_k\|_{E}=1$. Furthermore, we have that
\be\label{order}
\phi_k^2(x)\geq \f{\log k}{2M}+O(1)\ \   \text{for}\ \   |x|\leq \f{1}{\sqrt{k}}.
\ee

Now we can begin the proof of Theorem 1.2.\vspace{.2cm}

\noindent{\bf Proof of Theorem 1.2.} By $(H_2)$, we have $F(x,t\phi_k)\geq 0$ for all $t\geq 0$ and $x\in \mathbb{R}^{2m}$. This implies that
$$\int_{\mathbb{R}^{2m}} \f{F(x,t\phi_k)}{|x|^\beta}dx\geq
\int_{|x|\leq \f{1}{\sqrt{k}}} \f{F(x,t\phi_k)}{|x|^\beta}dx.$$
By (\ref{order}), we have, for $t$ and $k$ sufficiently large, there exists a constant $C_\phi>0$ such that
$$t\phi_k>C_\phi \ \ \text{for}\ \ |x|\leq \f{1}{\sqrt{k}}.$$
Since $(H_2)$ implies that, for $s>\f{C_\phi}{2}$,
$$\int_{\f{C_\phi}{2}}^s \f{\mu}{t}dt\leq \int_{\f{C_\phi}{2}}^s \f{f(x,t)}{F(x,t)}dt,$$
we have, if $t$ and $k$ sufficiently large, for $|x|\leq \f{1}{\sqrt{k}}$,
$$F(x,t\phi_k)\geq C t^\mu\phi_k^\mu.$$
Therefore,
\bna
J_{\epsilon}(t\phi_k)&=&\f{t^2}{2}
\int_{\mathbb{R}^{2m}}\le(|\nabla^m \phi_k|^2+\sum_{\gamma=0}^{m-1}a_\gamma (x)|\nabla^\gamma \phi_k|^2\ri)dx
-\int_{\mathbb{R}^{2m}}\f{F(x,t\phi_k)}{|x|^\beta}dx
-\epsilon t\int_{\mathbb{R}^{2m}}h\phi_kdx\\
&\leq& \f{t^2}{2}
\int_{\mathbb{R}^{2m}}\le(|\nabla^m \phi_k|^2+\sum_{\gamma=0}^{m-1}a_\gamma (x)|\nabla^\gamma \phi_k|^2\ri)dx
-\int_{|x|\leq \f{1}{\sqrt{k}}}\f{F(x,t\phi_k)}{|x|^\beta}dx
-\epsilon t\int_{\mathbb{R}^{2m}}h\phi_kdx\\
&\leq& \f{t^2}{2}
\int_{\mathbb{R}^{2m}}\le(|\nabla^m \phi_k|^2+\sum_{\gamma=0}^{m-1}a_\gamma (x)|\nabla^\gamma \phi_k|^2\ri)dx
-Ct^\mu\int_{|x|\leq \f{1}{\sqrt{k}}}\f{\phi_k^\mu}{|x|^\beta}dx
-\epsilon t\int_{\mathbb{R}^{2m}}h\phi_kdx
\ena
Since $\mu>2$, we get
\be\label{functionalinfty}
\lim_{t\ra +\infty}J_{\epsilon}(t\phi_k)=-\infty.
\ee

Suppose (\ref{minmaxlevel}) is not correct. Then we have, for all $k$ and $\epsilon>0$,
\be\label{antiminmax}
\max_{t\geq 0}J_{\epsilon}(t\phi_k)\geq \le(1-\f{\beta}{2m}\ri)\f{(4\pi)^m m!}{2\alpha_0}.
\ee
(\ref{functionalinfty}) and (\ref{antiminmax}) imply that, for any fixed $k$, there exists $t_k>0$ such that
$$J_{\epsilon}(t_k\phi_k)=\max_{t\geq 0}J_{\epsilon}(t\phi_k).$$
It follows that $\f{d}{dt}J_{\epsilon}(t\phi_k)=0$ at $t=t_k$, or equivalently,
\be\label{tk}
t_k^2=t_k^2\|\phi_k\|_E^2=\int_{\mathbb{R}^{2m}}
\f{t_k\phi_kf(x,t_k\phi_k)}{|x|^\beta}dx
+\epsilon t_k\int_{\mathbb{R}^{2m}}h\phi_kdx.
\ee

Now we claim that $\{t_k\}$ is a bounded sequence and its upper bound is independent of $\ep$. Suppose not. $(H_5)$ implies that, for any $\sigma>0$, there exists $R_\sigma>0$ such that, for all $s\geq R_\sigma$, it holds that
$$sf(x,s)\geq \sigma e^{\alpha_0 s^2}.$$
Then by (\ref{order}) and (\ref{tk}), we have, for sufficiently large $k$,
\bea\label{tkinequality}
t_k^2&\geq& \sigma\int_{|x|\leq\f{1}{\sqrt{k}}}\f{e^{\alpha_0t_k^2\phi_k^2}}{|x|^\beta}dx-\epsilon t_k\|h\|_{E^*}\|\phi_k\|_E\nonumber\\
&\geq& \sigma\int_{|x|\leq\f{1}{\sqrt{k}}}\f{e^{\alpha_0t_k^2\le(\f{\log k}{2M}+O(1)\ri)}}{|x|^\beta}dx-\epsilon t_k\|h\|_{E^*}\|\phi_k\|_E\nonumber\\
&=&\f{\sigma\omega_{2m-1}}{(2m-\beta)k^{m-\f{\beta}{2}}
}e^{\alpha_0t_k^2\le(\f{\log k}{2M}+O(1)\ri)}
-\epsilon t_k\|h\|_{E^*}\|\phi_k\|_E,
\eea
(\ref{tkinequality}) is equivalent to
\bna
1&\geq& \f{\sigma\omega_{2m-1}}{(2m-\beta)}k^{\f{\alpha_0t_k^2}{2M}+o(1)
-m+\f{\beta}{2}-\f{2\log t_k}{\log k}}-\ep\|h\|_{E^*}\|\phi_k\|_E k^{\f{\log t_k}{\log k}}\\
&\geq& \f{\sigma\omega_{2m-1}}{(2m-\beta)}k^{\f{\alpha_0t_k^2}{2M}+o(1)
-m+\f{\beta}{2}-\f{2\log t_k}{\log k}}-\|h\|_{E^*}\|\phi_k\|_E k^{\f{\log t_k}{\log k}}.
\ena
Let $k\ra \infty$, we get a contradiction because the right hand side of the inequality tends to $+\infty$. Thus the claim is proved.

From(\ref{antiminmax}), we get that
\be\label{tkminmax}
t_k^2\geq \le(1-\f{\beta}{2m}\ri)\f{(4\pi)^m m!}{\alpha_0}+2\epsilon t_k\int_{\mathbb{R}^{2m}}h\phi_kdx.
\ee
Here we have used the fact that $\|\phi_k\|_E=1$ and $F(x,s)\geq 0$.
Since
$$\le|\epsilon t_k\int_{\mathbb{R}^{2m}}h\phi_kdx\ri|\leq \epsilon t_k\|h\|_{E^*}\|\phi_k\|_E=\epsilon t_k\|h\|_{E^*},$$
$t_k$ is bounded and $\epsilon$ can be arbitrarily small, (\ref{tkminmax}) implies that
$$t_k^2\geq \le(1-\f{\beta}{2m}\ri)\f{(4\pi)^m m!}{\alpha_0}.$$
If
$$\lim_{k\ra\infty}t_k^2> \le(1-\f{\beta}{2m}\ri)\f{(4\pi)^m m!}{\alpha_0},$$
we get, for sufficiently large $k$,
$$\alpha_0 t_k^2\f{\log k}{2M}>\le(m-\f{\beta}{2}\ri)\log k.$$
This is a contradiction with the fact that $\{t_k\}$ is a bounded sequence because the right hand side of (\ref{tkinequality}) tends to $+\infty$ as $k\ra +\infty$. Thus we have
$$\lim_{k\ra\infty}t_k^2=\le(1-\f{\beta}{2m}\ri)\f{(4\pi)^m m!}{\alpha_0}.$$
Let $k\ra \infty$ and $\epsilon\ra 0$ in (\ref{tkinequality}), we obtain
$$\le(1-\f{\beta}{2m}\ri)\f{(4\pi)^m m!}{\alpha_0}\geq \f{\sigma\omega_{2m-1}}{(2m-\beta)}.$$
This is a contradiction because $\sigma$ can be chosen arbitrarily large. Thus Theorem 1.2 is proved under assumption $(H_5)$.

If $f(x,s)$ satisfies $(H_5)'$ instead of $(H_5)$. We can choose a bounded sequence of functions $\{u_k\}\subset E$ such that
$$\int_{\mathbb{R}^{2m}}\f{|u_k|^p}{|x|^\beta}dx=1\quad{\rm and}\quad \|u_k\|_E\ra S_p.$$
Then by Lemma 2.1, we can assume that there exists a function $u_p$ such that
   $$\begin{array}{lll}
   u_k\rightharpoonup u_p &{\rm in}\quad E,\\[1.5ex]
   u_k\ra u_p &{\rm in}\quad L^q(\mathbb{R}^{2m})\quad{\rm for\ \ all}\quad q\in [1,+\infty),\\[1.5ex]
   u_k(x)\ra u_p(x) &{\rm almost\ \ everywhere}.
   \end{array}$$
   These imply that
   $$\int_{\mathbb{R}^{2m}}\f{|u_k|^p}{|x|^\beta}dx\ra \int_{\mathbb{R}^{2m}}\f{|u_p|^p}{|x|^\beta}dx=1.$$
   On the other hand, we have
   $$\|u_p\|_E\leq \liminf_{k\ra\infty}\|u_k\|_E=S_p .$$
   Thus we get $\|u_p\|_E=S_p.$ Define a function $M_\ep(t): [0, +\infty)\ra \mathbb{R}$ by
    $$M_\ep(t):=\f{t^2}{2}\int_{\mathbb{R}^{2m}}\le(|\nabla^m u_p|^2+\sum_{\gamma=0}^{m-1}a_\gamma (x)|\nabla^\gamma u_p|^2\ri)dx
-\int_{\mathbb{R}^{2m}}\f{F(x,tu_p)}{|x|^\beta}dx
-t\epsilon\int_{\mathbb{R}^{2m}}hu_pdx$$
   By $(H_2)$, $(H_5)'$ and $\int_{\mathbb{R}^{2m}}\f{|u_p|^p}{|x|^\beta}dx=1$, we have
   \bna
   M_\ep(t)
   &\leq&\f{t^2}{2}\int_{\mathbb{R}^{2m}}\le(|\nabla^m u_p|^2+\sum_{i=0}^{m-1}a_i (x)|\nabla^i u_p|^2\ri)dx
   -C_p\f{t^p}{p}\int_{\mathbb{R}^{2m}}\f{|u_p|^p}{|x|^\beta}dx+\ep t \|h\|_{E^*}\|u_p\|\\
   &=&\f{t^2}{2}S_p^2-C_p\f{t^p}{p}+\ep t S_p\|h\|_{E^*}\\
   &\leq&\f{(p-2)}{2p}\f{S_p^{2p/(p-2)}}{C_p^{2/(p-2)}}+\ep t_0 S_p \|h\|_{E^*},
   \ena
where $t_0$ is a constant which belongs to $[0,+\infty)$ and is independent of the choice of $\ep$. By choosing $\ep$ small enough, we get the desired results from the definitions of $C_p$ and $S_p$ in $(H_5)'$ immediately.$\hfill\Box$\\

\noindent{\bf The sharpness of the Adams inequality.} 

Define $\varphi_k=\f{\tilde{\phi}_k}{\|\tilde{\phi}_k\|_{\tilde{E}}}.$
Obviously, we have
\bea\label{sharpness}
\sup_{u\in W^{m,2}(\mathbb{R}^{2m}),\|u\|_{\tilde{E}}\leq 1}
\int_{\mathbb{R}^{2m}}\f{e^{\alpha u^2}-1}{|x|^\beta}dx
&\geq& \int_{\mathbb{R}^{2m}}\f{e^{\alpha\varphi_k^2}-1}{|x|^\beta}dx\nonumber\\
&\geq&
\int_{|x|\leq\f{1}{\sqrt{k}}}\f{e^{\alpha\varphi_k^2}-1}{|x|^\beta}dx\nonumber\\
&\geq&
\int_{|x|\leq\f{1}{\sqrt{k}}}(Ck^{\f{\alpha}{2M}+\f{\beta}{2}}-k^{\f{\beta}{2}})dx\nonumber\\
&=&\f{\omega_{2m-1}(Ck^{\f{\alpha}{2M}+\f{\beta}{2}}-k^{\f{\beta}{2}})}{2m k^{m}}.
\eea
When $\alpha>\le(1-\f{\beta}{2m}\ri)\alpha(m,2m)$, substituting $M=\f{(4\pi)^m(m-1)!}{2}$, we get
$$\f{\alpha}{2M}+\f{\beta}{2}>m.$$
This implies that the right hand side of (\ref{sharpness}) tends to infinity as $k\ra \infty$. Thus the inequality (\ref{adams}) is sharp.

\section{Multiplicity result of the related elliptic equation}

To deal with equation (\ref{equation}), the main difference between our general case and the special case $n=2m=4$ is the function space $E$. As our proof of Lemma 2.1, the proofs of the following three lemmas are essentially the same as those in \cite{Yang}. The different definitions of $E$ do not cause difficulties and so we omit the proofs here.\vspace{.2cm}

\noindent{\bf Lemma 4.1.} {\it Assume $(A_1)$, $(A_2)$ and $(H_1)-(H_3)$. Then for any Palais-Smale sequence $\{u_k\}\subset E$ of $J_\ep$, i.e.,
$$J_\ep(u_k)\ra C, J'_\ep(u_k)\ra 0\ \ \text{as}\ \ k\ra \infty,$$
up to subsequence, there exists $u\in E$ such that
$$u_k\rightharpoonup u\ \ \text{in}\ \  E\ \ \text{and}\ \  u_k\ra u\ \ \text{in}\ \ L^q(\mathbb{R}^{2m})\ \ \text{for any}\ \ q\geq 1.$$
Furthermore, we have
$$\le\{\begin{array}{ll}
   \f{f(x,u_k)}{|x|^\beta}\ra\f{f(x,u)}{|x|^\beta}\ \ \text{in}\ \ L^1(\mathbb{R}^{2m}),\\[2.5ex]
   \f{F(x,u_k)}{|x|^\beta}\ra\f{F(x,u)}{|x|^\beta}\ \ \text{in}\ \ L^1(\mathbb{R}^{2m})
   \end{array}\ri.$$
and $u$ is a weak solution of (\ref{equation}).}\\

\noindent{\bf Lemma 4.2.} {\it Assume $(A_1)$, $(A_2)$, $(H_1)$, $(H_2)$ and $(H_4)$. Then there exists $\ep_2>0$ such that, for any $0<\ep<\ep_2$, there exists a Palais-Smale sequence
$\{u_k\}\subset E$ at level $C_0<0$ which converges strongly in $E$ to a minimum type solution $u_0$ of (\ref{equation}). Furthermore, we have $C_0\ra 0$ as $\ep\ra 0$. }\\

\noindent{\bf Lemma 4.3.} {\it Assume $(H_5)$ or $(H_5)'$, together with $(A_1)$, $(A_2)$ and  $(H_1)-(H_4)$. Then there exists $\ep_3>0$ such that, for any $0<\ep<\ep_3$, there exists a Palais-Smale sequence
$\{v_k\}\subset E$ at level $C_M>0$ which converges weakly in $E$ to a mountain-pass type solution $v_0$ of (\ref{equation}) with min-max level $C_M$.}\\

By Lemma 4.2 and 4.3, to prove Theorem 1.3, we only need to prove that $u_0$ and $v_0$ are distinct weak solutions. During the proof, we need the following singular version of Lions' inequality. This kind of inequality was first proved by Lions in \cite{Lions}.\vspace{.2cm}

\noindent{\bf Lemma 4.4.} {\it Let $\{w_k\}$ be a sequence in $E$. Suppose $\|w_k\|_E=1$ and $w_k\rightharpoonup w_0$ in $E$. Then, for any $0<p<\le(1-\f{\beta}{2m}\ri)\f{\alpha(m,2m)}{1-\|w_0\|_E^2}$, we have
$$\sup_k\int_{\mathbb{R}^{2m}}\f{e^{pw_k^2}-1}{|x|^\beta}dx<
\infty.$$}\\

\noindent{\bf Proof of Lemma 4.4.} If $w_0\equiv 0$, the lemma is a direct consequence of Theorem 1.1. Otherwise, by our assumptions on $w_k$, we have
\bea\label{lion4}
\|w_k-w_0\|_E^2
&=&1+\|w_0\|_E^2-2\int_{\mathbb{R}^{2m}}\le(\nabla^m w_k\nabla^m w_0+\sum_{\gamma=0}^{m-1}a_\gamma(x)\nabla^\gamma w_k\nabla^\gamma w_0\ri)dx\nonumber\\
&\ra& 1-\|w_0\|_E^2\ \ \text{as}\ \ k\ra \infty.
\eea
Using Young's inequality, we have, for any $\delta>0$,
\bea\label{lion3}
\int_{\mathbb{R}^{2m}}\f{e^{pw_k^2}-1}{|x|^\beta}dx
&\leq&\int_{\mathbb{R}^{2m}}
\f{e^{p\le((1+\delta)(w_k-w_0)^2+(1+\delta^{-1})w_0^2\ri)}-1}
{|x|^\beta}dx\nonumber\\
&\leq&\f{1}{\mu}\int_{\mathbb{R}^{2m}}
\f{e^{\mu p(1+\delta)(w_k-w_0)^2-1}}{|x|^\beta}dx
+\f{1}{\nu}\int_{\mathbb{R}^{2m}}
\f{e^{\nu p(1+\delta^{-1})w_0^2-1}}{|x|^\beta}dx\nonumber\\
&=:&W_1+W_2.
\eea
where $\mu>1$, $\nu>1$ and $\f{1}{\mu}+\f{1}{\nu}=1$.

We can choose $\delta$ sufficiently small and $\mu$ sufficiently close to $1$ such that
$$\mu p (1+\delta)(1-\|w_0\|_E^2)<\le(1-\f{\beta}{2m}\ri)\alpha(m,2m).$$
Then, by Theorem 1.1 and (\ref{lion4}), we have $W_1<C$ for some universal constant $C$.

To estimate $W_2$, we first claim that , for any $\alpha>0$ and $u\in E$, we have
$$\int_{\mathbb{R}^{2m}}\f{e^{\alpha u^2}-1}{|x|^\beta}dx<\infty.$$
In fact, since $E\subseteq W^{m,2}(\mathbb{R}^{2m})$, by the density of $C_0^\infty(\mathbb{R}^{2m})$ in $W^{m,2}(\mathbb{R}^{2m})$, we can choose some $u_0\in C_0^\infty(\mathbb{R}^{2m})$ such that
\be\label{lion1}
\|u-u_0\|_{W^{m,2}(\mathbb{R}^{2m})}^2
<\le(1-\f{\beta}{2m}\ri)\f{\alpha(m,2m)}{2\alpha}.
\ee
We get from (\ref{lion1}) that
\be\label{lion2}
\le\|\sqrt{\le(\f{2m}{2m-\beta}\ri)
\f{2\alpha}{\alpha(m,2m)}}(u-u_0)\ri\|_{W^{m,2}(\mathbb{R}^{2m})}
<1.
\ee
Assume $R_u$ and $C_u$ are positive constants such that supp$u_0\subset B_{R_u}$ and $|u_0|\leq C_u$. Then we have
\bna
\int_{\mathbb{R}^{2m}}\f{e^{\alpha u^2}-1}{|x|^\beta}dx
&\leq& \int_{\mathbb{R}^{2m}}\f{e^{2\alpha (u-u_0)^2+2\alpha u_0^2}-1}{|x|^\beta}dx\\
&=&\int_{\mathbb{R}^{2m}}\f{e^{2\alpha (u-u_0)^2+2\alpha u_0^2}-e^{2\alpha u_0^2}+e^{2\alpha u_0^2}-1}{|x|^\beta}dx\\
&\leq&e^{2\alpha C_u^2}\int_{\mathbb{R}^{2m}}\f{e^{2\alpha (u-u_0)^2}-1}{|x|^\beta}dx+\int_{B_{R_u}}\f{e^{2\alpha u_0^2}-1}{|x|^\beta}dx\\
&\leq&e^{2\alpha C_u^2}\int_{\mathbb{R}^{2m}}\f{e^{2\alpha (u-u_0)^2}-1}{|x|^\beta}dx+(e^{2\alpha C_u^2}-1)\int_{B_{R_u}}\f{1}{|x|^\beta}dx\\
&=&e^{2\alpha C_u^2}\int_{\mathbb{R}^{2m}}\f{e^{2\alpha (u-u_0)^2}-1}{|x|^\beta}dx+\f{\omega_{2m-1}(e^{2\alpha C_u^2}-1)}{2m-\beta}R_u^{2m-\beta}\\
\ena
By (\ref{adams}) with $\tau_\gamma=1$ for $0\leq \gamma\leq m$ and (\ref{lion2}), we have
$$\int_{\mathbb{R}^{2m}}\f{e^{2\alpha (u-u_0)^2}-1}{|x|^\beta}dx\leq C $$
for some universal constant $C$. Thus we get our claim proved and it is easy to see that $W_2<C$ follows from the claim immediately.
$\hfill\Box$\\

\noindent{\bf Proof of Theorem 1.3.} Suppose that $u_0$ and $v_0$ are the minimum and mountain-pass type solutions of (\ref{equation}) respectively. By Lemma 4.2 and 4.3, we have that, for $\ep$ small enough, there are two
 sequences
  $\{u_k\}$ and $\{v_k\}$ in $E$
  such that
$$u_k\ra u_0\quad {\rm and} \quad v_k\rightharpoonup
v_0,$$
$$J_{\epsilon}(u_k)\ra C_0<0\quad {\rm and} \quad J_{\epsilon}(v_k)
\ra C_M>0,$$
$$J_{\epsilon} '(u_k)u_k\ra 0\quad {\rm and} \quad
J_{\epsilon} '(v_k)v_k\ra 0.$$
Theorem 1.2 tells us that
$0<C_M<C_0+\le(1-\f{\beta}{2m}\ri)\f{(4\pi)^m m!}{2\alpha_0}.$ We will show a contradiction under the
assumption $u_0=v_0$.

Let
$$w_k=\f{v_k}{\|v_k\|_E}\quad {\rm and}
\quad w_0=\f{u_0}{\lim\limits_{k\ra\infty} \|v_k\|_E}.$$ We have
$\|w_k\|=1$ and $w_k\rightharpoonup w_0$ in $E$. In particular
$\|w_0\|_E\leq 1$.
To proceed, we distinguish two cases.\\

\noindent{\bf Case 1.} $\|w_0\|_E=1$.

In this case, we have
$$\lim\limits_{k\ra \infty}\|v_k\|_E
=\|u_0\|_E.$$
Therefore, $v_k\ra u_0$ in $E$. Lemma 4.1 tells us that $$\f{F(v_k)}{|x|^\beta}\ra \f{F(u_0)}{|x|^\beta}
\quad {\rm in} \quad L^1(\mathbb{R}^{2m}).$$ 
Then we have
$$J_{\epsilon}(v_k)\ra J_{\epsilon}(u_0)=C_0,$$
which is a contradiction with our assumption.\\

\noindent{\bf Case 2.} $\|w_0\|_E<1$.

Since $0<C_M<C_0+\le(1-\f{\beta}{2m}\ri)\f{(4\pi)^m m!}{2\alpha_0}=J_{\ep}(u_0)+\le(1-\f{\beta}{2m}\ri)\f{(4\pi)^m m!}{2\alpha_0}$, we can
   choose some $q>1$ sufficiently close to $1$ and
   $\delta>0$ such that
$$q\alpha_0 \|v_k\|_E^{2}\leq
\le(1-\f{\beta}{2m}\ri)\f{(4\pi)^m m!\|v_k\|_E^2}{2(C_M-J_{\ep}(u_0))}-\delta.$$ 
Since $v_k\rightharpoonup u_0$ in $E$ and $\f{F(x,v_k)}{|x|^\beta}\ra
\f{F(x,u_0)}{|x|^\beta}$ in $L^1(\mathbb{R}^{2m})$, we have
\bna
\lim_{k\ra \infty}\|v_k\|_E^2(1-\|w_0\|_E^2)=\lim_{k\ra
\infty}\|v_k\|_E^2-\|u_0\|_E^2=2(C_M-J_{\ep}(u_0)).
\ena
Then we get, for $k$
sufficiently large,
\be\label{vkw0}
q\alpha_0 \|v_k\|_E^{2}<
\le(1-\f{\beta}{2m}\ri)\f{(4\pi)^m m!}{1-\|w_0\|_E^2}.
\ee
Suppose $\alpha>0$, $p>1$ and $p'>p$, using L'Hospital's rule, we have that there exists a positive constant $C_\alpha$ which only depends on $\alpha$, such that for all $s>0$,
$$(e^{\alpha s^2}-1)^p\leq C_\alpha (e^{\alpha p's^2}-1).$$
In fact this is a result proved by the first author in \cite{Zhao1}.
Then by $(H_1)$ and H\"{o}lder's inequality, we have
\bna
\int_{\mathbb{R}^{2m}}\f{|f(x,v_k)|^p}{|x|^{\beta}}dx
&\leq&
C\int_{\mathbb{R}^{2m}}\f{|v_k|^p+ |v_k|^{p\theta}(e^{\alpha_0v_k^2}-1)^p}{|x|^{p\beta}}dx\\
&\leq&
C\int_{\mathbb{R}^{2m}}\f{|v_k|^p}{|x|^{p\beta}}dx
+C\int_{\mathbb{R}^{2m}}\f{|v_k|^{pp_1\theta}}{|x|^{p\beta}}dx
\int_{\mathbb{R}^{2m}}\f{(e^{\alpha_0v_k^2}-1)^{pp_2}}{|x|^{p\beta}}dx\\
&\leq&
C\int_{\mathbb{R}^{2m}}\f{|v_k|^p}{|x|^{p\beta}}dx
+C\int_{\mathbb{R}^{2m}}\f{|v_k|^{pp_1\theta}}{|x|^{p\beta}}dx
\int_{\mathbb{R}^{2m}}\f{e^{\alpha_0pp_2'v_k^2}-1}{|x|^{p\beta}}dx,
\ena
where $\f{1}{p_1}+\f{1}{p_2}=1$ and $p_2<p_2'$. Since $0\leq \beta<2m$, (\ref{vkw0}) and the continuous embedding $E\hookrightarrow L^p(\mathbb{R}^{2m})$ for any $p\geq 1$  imply that there exists some $p: 1<p<q$ such that $\f{f(x,v_k)}{|x|^\beta}$ is bounded in $L^{p}(\mathbb{R}^{2m})$. It follows that
\bna \left|\int_{\mathbb{R}^{2m}} \f{f(x,v_k)(v_k-u_0)}
{|x|^\beta}dx\right|
\leq C\|v_k-u_0\|_{L^{\f{p}{1-p}}}\ra 0. \ena
From this convergence and
$J_\ve '(v_k)(v_k-u_0)\ra 0$,
 we get
$$\int_{\mathbb{R}^{2m}} \le(\nabla^m v_k \nabla^m(v_k-u_0)+\sum_{\gamma=0}^{m-1}a_\gamma(x) \nabla^\gamma v_k\nabla^\gamma(v_k-u_0)\ri)dx \ra 0.$$ Moreover,
since $v_k\rightharpoonup u_0$, we have
$$\int_{\mathbb{R}^{2m}} \le(\nabla^m u_0 \nabla^m(v_k-u_0)+\sum_{\gamma=0}^{m-1}a_\gamma(x) \nabla^\gamma u_0\nabla^\gamma(v_k-u_0)\ri)dx \ra 0.$$
These two limitations tell us that $v_k\ra u_0$ in $E$. From the continuity of the functional $J_{\ep}$, we get $J_\ep(v_k)\ra J_\ep(u_0)=C_0$, which is still a contradiction and the proof is
finished.
$\hfill\Box$\\

{\bf Acknowledgements.} The first author is partially supported by NSFC (11001268, 11071020), Beijing Natural Science Foundation No.2112023 and PCSIRT.


\begin{thebibliography}{00}

\bibitem{AdaTan} S. Adachi, K. Tanaka, Trudinger type inequalities in $\mathbb{R}^N$ and their best exponents, Proc. Amer. Math. Soc., 128 (1999) 2051-2057.

\bibitem{Adams} D.R. Adams, A sharp inequality of J. Moser for higer order derivatives, Ann. of Math., 128 (1998) 385-398.

\bibitem{Adimurthi} Adimurthi, Existence of positive
solutions of the semilinear Dirichlet problem with critical growth
for the $n$-Laplacian, Ann. Scuola Norm. Sup. Pisa Cl. Sci., 17
(1990) 393-413.

\bibitem{AdiYad} Adimurthi, S.L. Yadava, Multiplicity results for semilinear elliptic equations in a bounded domain of $\mathbb{R}^2$ involving critical exponents, Ann. Scuola Norm. Sup. Pisa Cl. Sci., 17 (1990) 481-504.

\bibitem{AY} Adimurthi, Y.Y. Yang, An interpolation
of Hardy inequality and Trudinger-Moser inequality in $\mathbb{R}^N$
and its applications, International Mathematics Research Notices 13
(2010) 2394-2426.

\bibitem{Alb} A. Alberico, Moser type inequalities for higer-order derivatives in Lorentz spaces, Potential Anal., 28 (2008) 389-400.

\bibitem{Cao} D.M. Cao, Nontrivial solution of semilinear elliptic equation with
critical exponent in $\mathbb{R}^2$, Comm. Partial Differential
Equations, 17 (1992) 407-435.

\bibitem{doo} J.M. do \'O, Semilinear Dirichlet problems for the $N$-Laplacian in $\mathbb{R}^N$ with nonlinearities in critical growth range, Differential Integral Equations, 9 (1996) 967-979.

\bibitem{doo3} J.M. do \'O, E. Medeiros, U. Severo,
On a quasilinear nonhomogeneous elliptic equation with critical
growth in $\mathbb{R}^N$, J. Differential Equations 246 (2009)
1363-1386.

\bibitem{FDR} D.G. de Figueiredo, J.M. do \'O, B. Ruf, On an inequality by N. Trudinger and J. Moser and related elliptic equations, Comm. Pure Appl. Math., 25 (2002) 135-152.

\bibitem{FMR} D.G. de Figueiredo, O.H. Miyagaki, B. Ruf, Elliptic equations in $\mathbb{R}^2$ with nonlinearities in the critical growth range, Calc. Var. Partial Differential Equations, 3 (1995) 139-153.

\bibitem{Fon} L. Fontana, Sharp borderline Sobolev inequalities on compact Riemannian manifolds, Comm. Math. Helv., 68 (1993) 415-454.

\bibitem{GazGruSqu} F. Gazzola, H. Grunau, M. Squassina, Existence and nonexistence results for critical growth biharmonic elliptic equations, Calc. Var. Partial Differential Equations, 18 (2003) 117-143.

\bibitem{Kav} S. Kesavan, Symmetrization and applications, Series in Analysis, 3. World Scientific Publishing Co. Pte. Ltd., Hackensack, 2006.

\bibitem{LazSch} M. Lazzo, P.G. Schmidt, Oscillatory radial solutions for subcritical biharmonic equations, J. Differential Equations, 247 (2009) 1479-1504.

\bibitem{LL1} N. Lam, G.Z. Lu, Existence and multiplicity of solutions to equations of $N$-Laplacian type with critical exponential growth in $\mathbb{R}^N$, J. Funct. Anal., 262 (2012), 1132-1165.

\bibitem{LL3} N. Lam, G.Z. Lu, Existence of nontrivial solutions to polyharmonic equations with subcritical and critical exponential growth, Discrete Contin. Dyn. Syst., 32 (2012) 2187-2205.

\bibitem{LL4} N. Lam, G.Z. Lu, Sharp Adams type inequalities in Sobolev spaces $W^{m,\f{n}{m}}(\mathbb{R}^n)$ for arbitrary integer $m$, J. Differential Equations, 253 (2012) 1143-1171.

\bibitem{LL2} N. Lam, G.Z. Lu, Sharp singular Adams inequalities in high order Sobolev spaces, preprint.

\bibitem{LR} Y.X. Li, B. Ruf, A sharp Trudinger-Moser type inequality for unbounded domains in $\mathbb{R}^n$, Indiana Univ. Math. J., 57 (2008) 451-480.

\bibitem{Lions} P.L. Lions, The concentration-compactness principle in the calculus of variations, Part I, Revista Mathem\'{a}tica Iberoamericana, 1 (1985) 145-201.

\bibitem{Moser} J. Moser, A sharp form of an inequality by
N.Trudinger, Ind. Univ. Math. J. 20 (1971) 1077-1091.

\bibitem{Panda} R. Panda, On semilinear Neumann problems with critical growth for the $N$-Laplacian, Nonlinear Anal., 26 (1996) 1347-1366.

\bibitem{ReiWet} W. Reichel, W. Weth, Existence of solutions to nonlinear, subcritical higer order elliptic Dirichlet problems, J. Differential Equations, 248 (2010) 1866-1878.

\bibitem{Ruf} B. Ruf, A sharp Trudinger-Moser type inequality for unbounded domains in $\mathbb{R}^2$, J. Funct. Anal., 219 (2005) 340-367.

\bibitem{RufSani} B. Ruf, F. Sani, Sharp Adams-type inequalities in $\mathbb{R}^n$, to appear in Trans. Amer. Math. Soc.

\bibitem{Tarsi} C. Tarsi, Adams' inequality and limiting Sobolev embeddings into Zygmund spaces, Potential Anal., DOI: 10.1007/s11118-011-9259-4

\bibitem{Trudinger} N.S. Trudinger, On embeddings
into Orlicz spaces and some applications, J. Math. Mech. 17 (1967)
473-484.

\bibitem{Yang} Y.Y. Yang, Adams type inequalities and related elliptic partial differential equations in dimension four, J. Differential Equations 252 (2012) 2266-2295.

\bibitem{Yang1} Y.Y. Yang, Existence of positive solutions to quasi-linear elliptic equations with exponential growth in the whole
Euclidean space, J. Funct. Anal. 262 (2012) 1679-1704.

\bibitem{Yang2} Y.Y. Yang, Trudinger-Moser inequalities on complete noncompact Riemannian manifolds, arXiv:1112.0724v1.

\bibitem{YangZhao} Y.Y. Yang, L. Zhao, A class of
Adams-Fontana type inequalities and related functionals on
manifolds, Nonlinear Diff. Eqn. Appl. 17 (2009) 119-135.

\bibitem{Zhao1} L. Zhao, A multiplicity result for a singular and nonhomogeneous elliptic problem in $\mathbb{R}^n$, J. Partial Diff. Eq. 25 (2012) 90-102.

\bibitem{Zhao} L. Zhao, Exponential problem on a compact
 Riemannian manifold without boundary, Nonlinear Analysis 75 (2012) 433-443.

\end{thebibliography}
\end{document}